\documentstyle[amssymb]{article}

\newtheorem{theorem}{Theorem}

\newtheorem{definition}[theorem]{Definition}
\newtheorem{example}[theorem]{Example}

\def\QED{\quad\blackslug\lower 8.5pt\null}

\begin{document}

\begin{center}
{\Large \bf SINGULAR POINTS } 
\vspace*{2mm}

{\Large \bf OF LIGHTLIKE  HYPERSURFACES}

\vspace*{2mm}

{\Large \bf OF THE DE SITTER SPACE}

\vspace*{3mm}
{\large M. A. Akivis and  V. V. Goldberg}\footnote{This research was partially supported by 
the Volkswagen-Stiftung (RiP-program at Oberwolfach).}
\end{center}

\vspace*{5mm}

\begin{abstract}
The authors study singular points  of lightlike hypersurfaces of 
the de Sitter space $S^{n+1}_1$ and  the geometry of hypersurfaces 
and use them for construction of an invariant normalization and an 
invariant affine connection of lightlike hypersurfaces.
\end{abstract}

\setcounter{equation}{0}

\setcounter{section}{-1}

\section{Introduction}  

It is well-known that the pseudo-Riemannian manifolds $(M, g)$ 
of Lorentzian signature play a special role in geometry and 
physics, and that they are models of spacetime of 
general relativity. At the tangent space $T_x$ of an arbitrary 
point $x$ of such a manifold, one can invariantly define a real 
isotropic cone $C_x$. From the point of view of physics, 
this cone is the light cone: trajectories of light impulses 
emanating from the point $x$ are tangent to this cone. 

Hypersurfaces  of a Lorentzian manifold $(M, g)$ can be of 
three types: spacelike, timelike, and lightlike 
(see, for example,  \cite{ON83} or \cite{AG96}). 
For definiteness, we will assume that $\dim M = 
n + 1$ and $\mbox{{\rm sign}} \; g = (n, 1)$.

The tangent hyperplane to a spacelike hypersurface $U^n$ 
at any point does not have real 
common points with  the light cone $C_x$. This implies that 
on $U^n$ a proper Riemannian metric is induced. 
The tangent hyperplane to a timelike  hypersurface $U^n$ 
at any point  intersects  the light cone $C_x$ along 
an $(n-1)$-dimensional cone. 
This implies that on $U^n$ a pseudo-Riemannian metric 
of Lorentzian signature $(n - 1, 1)$ is induced. Finally, the 
tangent hyperplane to a lightlike hypersurface $U^n$ 
at any point is tangent to  the light cones $C_x$. This implies that on $U^n$ a degenerate Riemannian metric 
signature $(n - 1, 0)$ is induced. 

On spacelike and timelike hypersurfaces of a manifold 
of Lorentzian signature, an invariant normalization and 
an affine Levi-Civita  connection are induced by 
a first-order neighborhood while on lightlike hypersurfaces 
one should use differential neighborhoods of higher order 
to construct an invariant normalization and 
an affine  connection,

From the point of view of physics lightlike hypersurfaces 
are of great importance since they are models of different 
types of horizons studied in general relativity: 
event horizons, Cauchy's horizons, Kruskal's horizons 
(see \cite{Ch83} and \cite{MTW73}). This is the reason that 
the study of geometric structure of lightlike hypersurfaces  
is of interest.

In the current paper we consider lightlike hypersurfaces 
in the de Sitter space (a pseudo-Riemannian space 
of Lorentzian signature and constant positive curvature), 
study their  geometric structure, and prove that there are 
singular points and singular submanifolds on them. 

The de Sitter space $S_1^{n+1}$ admits a realization on the 
exterior of an $n$-dimensional oval hyperquadric $Q^n$ of 
a projective space $P^{n+1}$. Thus the de Sitter space 
is isometric to a pseudoelliptic space, 
$S_1^{n+1} \sim \mbox{{\rm ext}} \; Q^n$. 
Since the interior of the 
 hyperquadric $Q^n$ is isometric to the hyperbolic geometry 
of the Lobachevsky space $H^{n+1}$, $H^{n+1} \sim  
\mbox{{\rm int}} \; Q^n$ and the geometry of $Q^n$ itself is 
equivalent to that of an  $n$-dimensional conformal 
space $C^n$, $C^n \sim Q^n$, the groups of motions 
of these three spaces are isomorphic to each other 
and are isomorphic to the group ${\bf SO} (n + 2, 1)$ 
of rotations of a pseudo-Euclidean space $R^{n+2}_1$ 
of Lorentzian signature. This allows us to apply the apparatus 
developed in the book \cite{AG96} for the conformal space $C^n$ 
to the study of the de Sitter space. 

Note also that the geometry of lightlike hypersurfaces 
on pseudo-Riemannian manifolds of different signatures 
was the subject of  many journal papers and even two 
books \cite{DB96} and \cite{Ku96}. 
However, the geometry of lightlike 
hypersurfaces in the de Sitter space was not studied 
in spite of the fact that this geometry has many interesting 
geometric features.

In the present paper we study the geometry of the 
de Sitter space $S^{n+1}_1$ using its connection with the 
geometry of the conformal space. We prove that the geometry of 
lightlike hypersurfaces of the space $S^{n+1}_1$ is directly 
connected with the geometry of hypersurfaces of the conformal 
space $C^n$. The 
latter  was studied in detail in the papers of the first 
author (see \cite{A52},  \cite{A65}) 
and also in the book \cite{AG96}. This simplifies the 
study of lightlike hypersurfaces of the de Sitter space 
$S^{n+1}_1$ and makes possible to apply for their consideration 
the apparatus constructed in the conformal theory. 

In Section 1 we study the geometry of the 
de Sitter space and its connection with the geometry of the 
conformal space. Next we study lightlike hypersurfaces 
$U^n$ in the space $S^{n+1}_1$, investigate their structure, 
and prove that such a hypersurface is tangentially degenerate of 
rank $r \leq n - 1$. Its rectilinear or plane 
generators form an isotropic fibre bundle on $U^n$.

In Sections 2--5 we investigate lightlike hypersurfaces $U^n$ of 
maximal rank, and for their study we use the relationship between 
the geometry of such hypersurfaces and the geometry of 
hypersurfaces of the conformal space. For a lightlike 
hypersurface, we construct the fundamental quadratic forms and 
connections determined  by a normalization of a hypersurface by 
means of a distribution (the screen distribution) which is 
complementary to the isotropic distribution. The screen 
distribution plays an important role in the book \cite{DB96} 
since it defines a connection on a lightlike hypersurface $U^n$, 
and it appears to be important for applications. We prove that 
the screen distribution on a lightlike hypersurface  can be 
constructed invariantly by means of quantities from a third-order 
differential neighborhood, that is, such a distribution is 
intrinsically connected with the geometry of a hypersurface. 

In Section 5 we study singular points of a  lightlike 
hypersurface in the de Sitter space $S^{n+1}_1$, classify them, 
and describe the structure of hypersurfaces carrying singular 
points of different types. Moreover, we establish the connection 
of this classification with that of canal hypersurfaces of the 
conformal space.

The principal method of our investigation is the method of 
moving frames and exterior differential forms in the form 
in which it is presented in the books \cite{AG93} and 
\cite{AG96}. All functions considered in the paper are assumed to 
be real and 
differentiable, and all manifolds are assumed to be smooth 
with the possible exception of some isolated singular points and 
singular submanifolds.

\section{The de Sitter Space}

\setcounter{equation}{0}
 
{\bf 1.} 
In a projective space $P^{n+1}$ of dimension $n+1$ we consider 
an oval hyperquadric $Q^n$. Let $x$ be a point of the space 
$P^{n+1}$ with projective coordinates $(x^0, x^1, \ldots , 
x^{n+1})$.  The hyperquadric $Q^n$ is determined by the equations 
\begin{equation}\label{eq:1}
(x, x) := g_{\xi\eta} x^\xi x^\eta = 0, \;\;\;\;\; 
\xi, \eta = 0, \ldots, n+1,
\end{equation}
whose left-hand side is a quadratic form $(x, x)$ 
of signature $(n+1, 1)$. The hyperquadric $Q^n$ 
divides the space $P^{n+1}$ into two parts, external and 
internal. Normalize the quadratic form $(x, x)$ 
in such a way that for the points of the external part 
the inequality $(x, x) > 0$ holds. This external domain is 
a model of the {\em de Sitter space} $S^{n+1}_1$ 
(see \cite{Z96}). We will identify the external domain 
of $Q^n$ with the space $S^{n+1}_1$. The hyperquadric $Q^n$ 
is the {\em absolute} of the space $S^{n+1}_1$.

On the hyperquadric $Q^n$ of  the space $P^{n+1}$ the geometry 
of a conformal space $C^n$ is realized. The bijective mapping 
$C^n \leftrightarrow Q^n$ is called the {\em Darboux mapping}, 
and the hyperquadric $Q^n$ itself is called the {\em Darboux 
hyperquadric}. 

Under the Darboux mapping to hyperspheres of the space $C^n$ 
there correspond cross-sections of the hyperquadric $Q^n$ 
by hyperplanes $\xi$. But to a hyperplane $\xi$ there corresponds 
a point $x$ that is polar-conjugate to $\xi$ with 
respect to $Q^n$ and lies outside of $Q^n$, that is, 
a point of the space $S^{n+1}_1$. Thus to hyperspheres of the 
space $C^n$ there correspond points  of the space $S^{n+1}_1$. 

Let $x$ be an arbitrary point of the space $S^{n+1}_1$. 
The tangent lines from the point $x$ to the hyperquadric $Q^n$ 
form a second-order cone $C_x$ with vertex at the point $x$. 
This cone is called the {\em isotropic cone}. 
For spacetime whose model 
is the space $S^{n+1}_1$ this cone is the light cone, and 
its generators are lines of propagation of light impulses 
whose source coincides with  the point $x$. 

The cone $C_x$ separates all straight lines passing through 
 the point $x$ into spacelike (not having common points with 
the hyperquadric $Q^n$), timelike (intersecting $Q^n$ 
in two different points), and lightlike (tangent to $Q^n$). 
The  lightlike straight lines are generators of the cone $C_x$. 

To a spacelike straight line $l \subset S_1^{n+1}$ there 
corresponds an elliptic pencil of hyperspheres in the conformal 
space $C^n$. All hyperspheres of this pencil pass through a 
common $(n-2)$-sphere $S^{n-2}$ (the center of this pencil). 
The sphere $S^{n-2}$ is the intersection of the hyperquadric 
$Q^n$ and the $(n-1)$-dimensional subspace of the space $P^{n+1}$ 
which is polar-conjugate to the line $l$ with respect to the 
hyperquadric $Q^n$. 

To a timelike straight line $l \subset 
S_1^{n+1}$ there corresponds a hyperbolic pencil of hyperspheres 
in the space $C^n$. 
Two arbitrary hyperspheres of this pencil do 
not have common points, and the pencil contains two hyperspheres 
of zero radius which correspond to the points of intersection of 
the straight line $l$ and the hyperquadric $Q^n$. 

Finally, to a lightlike straight line $l 
\subset S_1^{n+1}$ there corresponds a parabolic pencil of 
hyperspheres in the  space $C^n$ consisting of 
hyperspheres tangent one to another at a point that is a unique 
hypersphere of zero radius belonging to this pencil.

Hyperplanes of the space $S^{n+1}_1$ are also divided 
into three types. Spacelike hyperplanes do not have 
common points with  the hyperquadric $Q^n$; a timelike 
hyperplane intersects   $Q^n$ along a real hypersphere; 
and lightlike hyperplanes are tangent to $Q^n$. 
Subspaces of any dimension $r, \; 2 \leq r \leq n-1$, 
can be also classified in a similar manner.

Let us apply the method of moving frames 
to study some questions of differential geometry 
of the space $S^{n+1}_1$. With a point $x \in S^{n+1}_1$ 
we associate a family of projective frames $\{A_0, A_1, 
\ldots , A_{n+1}\}$. However, in order to apply formulas derived 
in the book \cite{AG96}, we will use the notations used 
in this book. Namely, we denote by $A_n$ the vertex 
of the moving frame which coincides 
 with the point $x$, $A_n = x$; we locate the vertices 
$A_0, A_i$, and $A_{n+1}$ at the hyperplane $\xi$ which is polar 
conjugate to the point $x$ with respect to 
 the hyperquadric $Q^n$, and we assume that the points $A_0$ and 
$A_{n+1}$ lie on the hypersphere $S^{n-1} = Q^n \cap \xi$, and 
the points $A_i$ are polar-conjugate to the straight line $A_0 
A_{n+1}$ with respect to $S^{n-1}$. Since $(x, x) > 0$, we can 
normalize the point $A_n$ by the condition $(A_n, A_n) = 1$. 
The points  $A_0$ and $A_{n+1}$ are not polar-conjugate  with 
respect to  the hyperquadric $Q^n$. Hence we can normalize them 
by the condition $(A_0, A_{n+1}) = - 1$. As a result, 
the matrix of scalar products of the frame elements 
has the form 
\begin{equation}\label{eq:2}
\renewcommand{\arraystretch}{1.3}
(A_\xi, A_\eta) = \pmatrix{0 & 0 & 0 & -1 \cr
0 & g_{ij} & 0 & 0 \cr
 0 & 0 & 1 & 0 \cr
-1 & 0 & 0& 0\cr}, \;\;\;\;\; i, j = 1, \ldots, n -1,
\renewcommand{\arraystretch}{1}
\end{equation}
and the quadratic form $(x, x)$ takes the form
\begin{equation}\label{eq:3}
(x, x) = g_{ij} x^i x^j + (x^n)^2 - 2 x^0 x^{n+1}.
\end{equation}
 The quadratic form $g_{ij} x^i x^j $ occurring in 
(3) is positive definite. 

The equations of infinitesimal displacement of the conformal 
frame $\{A_\xi\}, \linebreak \xi = 0, 1, \ldots, n+1$, 
we have constructed have the form
\begin{equation}\label{eq:4}
d A_\xi = \omega_\xi^\eta A_\eta, \;\;\;\; \xi, \eta 
= 0, 1, \ldots, n+1,
\end{equation}
where by (2), the 1-forms $\omega_\xi^\eta$  
  satisfy the following Pfaffian equations: 
\begin{equation}\label{eq:5}
\left\{
\renewcommand{\arraystretch}{1.3}
\begin{array}{ll}
\omega_0^{n+1} =  \omega^0_{n+1} = 0, &  \omega_0^0  + \omega_{n+1}^{n+1} = 0, \\ 
 \omega_i^{n+1} = g_{ij} \omega_0^j, & 
\omega_i^0 = g_{ij} \omega_{n+1}^j, \\ 
\omega_n^{n+1} - \omega_0^n   = 0, & \omega_n^0 - \omega_{n+1}^n = 0, \\
 g_{ij} \omega_n^j + \omega_i^n = 0 , &   \omega_n^n = 0, \\ 
 dg_{ij} = g_{jk} \omega_i^k + g_{ik} \omega_j^k. 
\end{array}
\renewcommand{\arraystretch}{1}
\right.
\end{equation}
These formulas are precisely the formulas derived in the book 
\cite{AG96} (see p. 32) for the conformal space $C^n$. 

 It follows from (4) that 
\begin{equation}\label{eq:6}
dA_n = \omega_n^0 A_0 + \omega_n^i A_i + \omega_n^{n+1} A_{n+1}.
\end{equation}

The differential $dA_n$ belong to the tangent space 
$T_x (S_1^{n+1})$, and the 1-forms $\omega_n^0, \omega_n^i$, 
and $\omega_n^{n+1}$ form a coframe of this space. The total 
number of these forms is $n + 1$, and this number coincides with 
the dimension of $T_x (S_1^{n+1})$. The scalar square of the 
differential $dA_n$ is the metric quadratic form $\widetilde{g}$ on 
the manifold $S_1^{n+1}$. By (2), this quadratic form $\widetilde{g}$ 
can be written as 
$$
\widetilde{g} = (dA_n, dA_n) = g_{ij} \omega_n^i \omega_n^j - 2 \omega_n^0 
\omega_n^{n+1}.
$$
Since the first term of this expression is a positive definite 
quadratic form, the form $\widetilde{g}$ is of Lorentzian signature 
$(n, 1)$. The coefficients of the form $\widetilde{g}$ 
produce the metric tensor of the space $S^{n+1}_1$ whose matrix is obtained from the 
matrix (2) by deleting the $n$th row and the $n$th column. 

The quadratic form $\widetilde{g}$  defines 
on $S^{n+1}_1$ a pseudo-Riemannian metric of signature $(n, 1)$. 
The isotropic cone defined in the space 
$T_x (S_1^{n+1})$ by the equation $\widetilde{g} = 0$ coincides with the cone $C_x$ 
that we defined earlier in the space $S^{n+1}_1$ geometrically. 

The 1-forms $\omega_\xi^\eta$ occurring in equations (4) 
  satisfy the structure equations of the space $C^n$:
\begin{equation}\label{eq:7}
d \omega_\xi^\eta = \omega_\xi^\zeta \wedge \omega_\zeta^\eta,
\end{equation}
which are obtained by taking exterior derivatives of 
equations (4) and which are conditions of complete 
integrability of (4). The forms $\omega_\xi^\eta$ 
are invariant forms of the fundamental group 
${\bf PO} (n+2, 1)$ of transformations of the 
spaces $H^{n+1}, C^n$, and $S^{n+1}_1$ which is locally 
isomorphic to the group ${\bf SO} (n+2, 1)$.

\section{Lightlike Hypersurfaces in the de Sitter Space}

A hypersurface $U^n$ in the de Sitter space $S^{n+1}_1$ 
is said to be {\em lightlike} if all its tangent hyperplanes are 
lightlike, that is, they are tangent to the hyperquadric $Q^n$ 
which is the absolute of the  space $S^{n+1}_1$.

Denote by $x$ an arbitrary point of the  hypersurface $U^n$, 
by $\eta$ the tangent hyperplane to $U^n$ at the point $x,  
\eta = T_x (U^n)$, and by $y$ the point of tangency of 
the hyperplane $\eta$ with the hyperquadric $Q^n$. 
Next, as in Section 1,  denote by $\xi$ the hyperplane 
which is polar-conjugate to the point $x$ 
with respect to the hyperquadric $Q^n$, and  associate 
with a point $x$ a family of projective  frames such that 
$x = A_n, y = A_0$, the points $A_i, i = 1, \ldots , 
n - 1$, belong to the intersection of the hyperplanes $\xi$ and 
$\eta, \; A_i \in \xi \cap \eta$, and the point $A_{n+1}$, as 
well as the point $A_0$, is the intersection point of the 
quadric $\xi \cap Q^n$ and the 
straight line that is polar-conjugate to the $(n-2)$-dimensional subspace spanned by the points $A_i$. 
In addition, we normalize the frame vertices 
in the same way as this was done in Section 1. Then the matrix of scalar products of the frame elements has the form (2), and the components of infinitesimal displacements of the moving frame 
satisfy the Pfaffian equations (5). 

Since the hyperplane $\eta$ is tangent to 
the  hypersurface $U^n$ at the point $x = A_n$ and does not 
contain the point $A_{n+1}$, the differential of 
 the point $x = A_n$ has the form
\begin{equation}\label{eq:8}
dA_n =\omega_n^0 A_0 + \omega_n^i A_i,
\end{equation}
 the following equation holds:
\begin{equation}\label{eq:9}
\omega_n^{n+1} = 0,
\end{equation}
and the forms $\omega_n^0$ and $\omega_n^i$ are basis forms of 
the  hypersurface $U^n$.
By relations (5), it follows from equation (8) that 
\begin{equation}\label{eq:10}
\omega_0^n = 0
\end{equation}
and 
\begin{equation}\label{eq:11}
dA_0 =\omega_0^0 A_0 + \omega_0^i A_i.
\end{equation}

Taking  exterior derivative of equation (9), we obtain 
\begin{equation}\label{eq:12}
\omega^i_n \wedge \omega_i^{n+1} = 0.
\end{equation}
Since the forms $\omega^i_n$ are linearly independent, 
by Cartan's lemma, we find from (12) that 
\begin{equation}\label{eq:13}
\omega_i^{n+1} = \nu_{ij} \omega_n^j, \;\; 
\nu_{ij} = \nu_{ji}.
\end{equation}
Applying an appropriate formula from (5), we find that
\begin{equation}\label{eq:14}
\omega_0^i = g^{ij}  \omega_j^{n+1} = 
g^{ik} \nu_{kj} \omega_n^j,
\end{equation}
where $(g^{ij})$ is the inverse matrix of the matrix $(g_{ij})$. 

Now formulas (8) and (11) imply that for 
$\omega_n^i = 0$, the point $A_n$ of 
the  hypersurface $U^n$ moves along the isotropic straight line 
$A_n A_0$, and hence $U^n$ is a ruled hypersurface. In what 
follows, we assume that the {\em entire}  straight line 
$A_n A_0$ belongs to the  hypersurface $U^n$. Thus the following 
theorem holds:

\begin{theorem}
A lightlike hypersurface $U^n$ is the image of the direct 
product $M^{n-1} \times l$ of a differentiable manifold 
$M^{n-1}$ and a projective line $l$ under the mapping $f: 
M^{n-1} \times l \rightarrow P^{n+1}$ into 
a projective space $P^{n+1}$: $U^n = f (M^{n-1} \times l)$ which sends 
the straight line $l$ to the straight line $A_n A_0 \in P^{n+1}$. 
\end{theorem}
Precisely this mapping is the subject of this paper.

In addition, formulas (8) and (11) show that at any point of 
a generator of  the  hypersurface $U^n$, its tangent hyperplane 
is fixed and coincides with the hyperplane $\eta$. Thus 
   $U^n$ is a {\em tangentially degenerate   hypersurface}. 

We recall that the {\em rank} of a tangentially degenerate   
hypersurface is the number of parameters on which the family of 
its tangent hyperplanes depends (see, for example, \cite{AG93}, 
p. 113). From relations (8) and (11) it follows that the 
tangent hyperplane $\eta$ of the hypersurface $U^n$ along 
its generator $A_n A_0$ is determined by this generator 
and the points $A_i$, $\eta = A_n \wedge A_0 \wedge A_1 
\wedge \ldots \wedge A_{n-1}$. The displacement of 
this hyperplane is determined by the differentials 
(8), (11), and
$$
dA_i = \omega_i^0 A_0 + \omega_i^j A_j + \omega_i^n A_n 
+ \omega_i^{n+1} A_{n+1}.
$$
But by (5), $\omega_i^n = - g_{ij} \omega_n^j$, and the forms 
$\omega_i^{n+1}$ are expressed according to formulas (13). 
From formulas (13) and (14) it follows that the rank of 
 a tangentially degenerate   hypersurface $U^n$ is 
determined by the rank of the matrix $(\nu_{ij})$ in 
terms of which the 1-forms $\omega_i^{n+1}$ and $\omega_0^i$ 
are expressed. But by (11) and (14) the dimension of 
the submanifold $V$ described by the point $A_0$ on the 
hyperquadric $Q^n$ is also equal to the rank of the matrix 
$(\nu_{ij})$. Thus we have proved 
the following result:

\begin{theorem} A lightlike hypersurface of the de Sitter 
space $S^{n+1}_1$ is a ruled tangentially degenerate 
hypersurface whose rank is equal to the dimension of 
the submanifold $V$ described by the point $A_0$ on the 
hyperquadric $Q^n$. 
\end{theorem}

Denote the rank of the tensor $\nu_{ij}$ and of  the hypersurface 
$U^n$ by $r$.  In this and next sections  we will assume that 
$r  = n - 1$. The case $r < n -1$ was considered by the authors 
in \cite{AG98b}. 

For  $r = n - 1$, the hypersurface $U^n$ carries an 
$(n-1)$-parameter 
family of isotropic rectilinear generators $l = A_n A_0$ along 
which the tangent hyperplane $T_x (U^n)$ is fixed. From the point 
of view of physics, the isotropic rectilinear generators of 
a lightlike hypersurface $U^n$ are trajectories of light 
impulses, and the hypersurface $U^n$ itself represents a 
{\em light flux} in spacetime.

Since $\mbox{{\rm rank}} \; (\nu_{ij}) = n - 1$,
the submanifold $V$ described by the point $A_0$ on the 
hyperquadric $Q^n$ has dimension $n-1$, that is, 
$V$ is a hypersurface. We denote it by $V^{n-1}$. The tangent 
subspace $T_{A_0} (V^{n-1})$ to $V^{n-1}$ is determined by 
the points $A_0, A_1, \ldots , A_{n-1}$. Since $(A_n, A_i) 
= 0$, this tangent subspace is polar-conjugate to the 
rectilinear generator $A_0 A_n$ of 
the lightlike hypersurface $U^n$.

The submanifold $V^{n-1}$ of the 
hyperquadric $Q^n$ is the image of a hypersurface of 
the conformal space $C^n$ under the Darboux mapping. 
 We will denote this hypersurface also by $V^{n-1}$. 
In the space $C^n$, 
the hypersurface  $V^{n-1}$ is defined by equation (10) which 
by (5) is equivalent to equation (9) defining a lightlike 
hypersurface $U^n$ in the space $S_1^{n+1}$. To 
the rectilinear generator $A_n A_0$ of the hypersurface $U^n$ 
there corresponds a parabolic pencil of hyperspheres $A_n + sA_0$ 
tangent to the hypersurface  $V^{n-1}$ (see \cite{AG96}, p. 40). 
Thus the following theorem is valid:

\begin{theorem} There exists a one-to-one correspondence 
between the set of hypersurfaces of the conformal space $C^n$ 
and the set of lightlike hypersurfaces of the maximal rank 
$r = n-1$ of the de Sitter 
space $S^{n+1}_1$. To  pencils of tangent hyperspheres 
of the hypersurface  $V^{n-1}$ there correspond isotropic 
rectilinear generators of the lightlike hypersurface $U^n$.
\end{theorem}

Note that for lightlike hypersurfaces of the four-dimensional 
Minkowski space $M^4$ the result similar to the result of Theorem 
2 was obtained  in \cite{K89}. 

\section{The Fundamental Forms and Connections on a 
Lightlike Hypersurface of the de Sitter Space}

The first fundamental form of a lightlike hypersurface 
$U^n$ of the space $S^{n+1}_1$ is a metric quadratic form. It 
is defined by the scalar square of the differential $dx$ 
of a point of this  hypersurface. Since we have $x = A_n$, 
by (8) and (2) this scalar square has the form 
\begin{equation}\label{eq:15}
(dA_n, dA_n) = g_{ij} \omega_n^i  \omega_n^j = g
\end{equation}
and is a positive semidefinite 
differential quadratic form of signature 
$(n-1, 0)$. It follows that the system of equations 
$ \omega_n^i = 0$ defines on the  hypersurface  $U^n$ 
a fibration of isotropic lines which, as we showed in Section 2, 
coincide with rectilinear generators of this hypersurface.

The second fundamental form of a lightlike hypersurface 
$U^n$ determines its deviation from the tangent hyperplane 
$\eta$. To find this quadratic form, we compute the part of 
the second differential of the point $A_n$ which 
does not belong to the tangent hyperplane $\eta = A_0 \wedge A_1
\wedge \ldots \wedge A_n$:
$$
d^2 A_n \equiv \omega_n^i \omega_i^{n+1} A_{n+1} \pmod{\eta}.
$$
This implies that the second fundamental form can be written as 
\begin{equation}\label{eq:16}
b = \omega_n^i  \omega_i^{n+1} = \nu_{ij} \omega_n^i  \omega_n^j,
\end{equation}
where we used expression (13) for the form $\omega_i^{n+1}$. 
Since we assumed that $\mbox{{\rm rank}} \;
(\nu_{ij}) = n - 1$, the rank of the quadratic form (16) 
as well as the rank of the form $g$ is equal to $n - 1$. 
The nullspace of this   quadratic form  (see \cite{ON83}, p. 53) 
is again determined by  the system of equations 
$ \omega_n^i = 0$  and coincides with the isotropic direction 
on the hypersurface $U^n$. The reduction of the rank of the 
quadratic form $b$ is connected with the tangential degeneracy 
of the hypersurface $U^n$. The latter was noted in Theorem 2.

On a hypersurface $V^{n-1}$ of the conformal space $C^n$ 
that corresponds to a lightlike hypersurface $U^n \subset 
S^{n+1}_1$, the quadratic forms (15) and (16) define 
the net of curvature lines, that is, an orthogonal 
and conjugate net.

To find the connection forms of the hypersurface $U^n$, 
we find exterior derivatives of its basis forms 
$\omega_n^0$ and $\omega_n^i$:
\begin{equation}\label{eq:17}
\renewcommand{\arraystretch}{1.3}
\left\{
\begin{array}{ll}
d \omega_n^0 =  \omega_n^0 \wedge \omega_0^0 
                  + \omega_n^i \wedge \omega_i^0, \\
d \omega_n^i =  \omega_n^0 \wedge \omega_0^i 
                 + \omega_n^j \wedge \omega_j^i.
\end{array}
\renewcommand{\arraystretch}{1}
\right.
\end{equation}
This implies that the matrix 1-form 
\begin{equation}\label{eq:18}
\omega = \pmatrix{\omega_0^0 & \omega_i^0 \cr
                   \omega_0^i & \omega_j^i \cr}
\end{equation} 
defines a torsion-free connection 
on the hypersurface $U^n$. To clarify the properties 
 of this connection, we find its curvature forms. 
Taking exterior derivatives of the forms 
(18) and applying equations (5), (7), (9), and (10), 
 we obtain 
\begin{equation}\label{eq:19}
\renewcommand{\arraystretch}{1.3}
\left\{
\begin{array}{lll}
\Omega_0^0 = d \omega_0^0 - \omega_0^i \wedge \omega_i^0 = 0, \\
\Omega_0^i = d \omega_0^i - \omega_0^0 \wedge \omega_0^i 
                   -  \omega_0^j \wedge \omega_j^i = 0, \\
\Omega_i^0 = d \omega_i^0 - \omega_i^0 \wedge \omega_0^0 
                   -  \omega_i^j \wedge \omega_j^0 
= -g_{ij}  \omega_n^j \wedge \omega_n^0, \\
\Omega_j^i = d \omega_j^i - \omega_j^0 \wedge \omega_0^i 
                   -  \omega_j^k \wedge \omega_k^i  
                   -  \omega_j^{n+1} \wedge \omega_{n+1}^i  
= -g_{jk}  \omega_n^k \wedge \omega_n^i. 
\end{array}
\renewcommand{\arraystretch}{1}
\right.
\end{equation}
In these formulas the forms $\omega_j^{n+1}$ and $\omega_0^i$ 
are expressed in terms of the basis forms $\omega_n^i$, and the 
forms $\omega_0^j, \omega_j^i$, and $\omega_i^0$ are 
fiber forms. If the principal parameters are fixed, then 
these fiber forms are invariant forms of the group $G$ 
of admissible transformations of frames associated with 
a point $x = A_n$ of the hypersurface $U^n$, and the connection 
 defined by the form (18) is a $G$-connection. 

To assign an affine connection on the hypersurface $U^n$, 
it is necessary to make a reduction of 
the family of frames in such a way that 
the forms $\omega_i^0$ become principal. Denote by 
$\delta$ the symbol of differentiation with respect to 
the fiber parameters, that is, for a fixed point $x = A_n$ 
of the hypersurface $U^n$, and by $\pi_\eta^\xi$ the values 
of the 1-forms $\omega_\eta^\xi$ for a fixed point $x = A_n$, 
that is, $\pi_\eta^\xi = \omega_\eta^\xi (\delta)$. Then 
we obtain
$$
\pi_n^0 = 0, \pi_n^i = 0, \pi_i^n = 0, \pi_i^{n+1} = 0.
$$
It follows 
\begin{equation}\label{eq:20}
\delta A_i = \pi_i^0 A_0 + \pi_i^j A_j.
\end{equation}

The points $A_0$ and $A_i$ determine the tangent subspace to the 
submanifold $V^{n-1}$ described by the point $A_0$ on the
 hyperquadric $Q^n$. If we fix an $(n-2)$-dimensional subspace 
$\zeta$ not containing the point $A_0$ in this tangent subspace 
and place the points $A_i$ into $\zeta$, then we obtain $\pi_i^0 = 0$. 
This means that the forms $\omega_i^0$ become principal, that is, 
\begin{equation}\label{eq:21}
\omega_i^0 = \mu_{ij} \omega_n^j + \mu_i \omega_n^0,
\end{equation}
and as a result, an affine connection arises 
on the hypersurface $U^n$. 

We will call the subspace $\zeta \subset T_{A_0} (V^{n-1})$ 
the {\em normalizing subspace} of the lightlike 
hypersurface $U^n$. We have proved the following result:

\begin{theorem} If in every tangent subspace $T_{A_0} (V^{n-1})$ 
of  the submanifold $V^{n-1}$ associated with a lightlike 
hypersurface $U^n, V^{n-1} \subset Q^n$, a normalizing 
 $(n-2)$-dimensional subspace $\zeta$  not containing the point 
$A_0$ is assigned, then there arises a torsion-free 
affine connection on $U^n$. 
\end{theorem}

The last statement follows the first two equations of (19).

The constructed above fibration of normalizing subspaces $\zeta$ 
defines a distribution $\Delta$ of $(n-1)$-dimensional 
elements on a lightlike hypersurface $U^n$. In fact, 
the point $x = A_n$ of the hypersurface $U^n$ along with the 
subspace $\zeta =  A_1 \wedge \ldots \wedge 
A_{n-1}$ define the $(n-1)$-dimensional subspace 
which is complementary to the straight line $A_n A_0$ and 
lies in the tangent subspace $\eta$ of the hypersurface $U^n$. 
Following the book \cite{DB96}, we will call this subspace 
the {\em screen}, and the distribution $\Delta$ the 
{\em screen distribution}. Since at the point $x$ the screen 
is determined by the subspace $A_n A_1 \ldots A_{n-1}$, 
the differential equations of the screen distribution 
has the form
\begin{equation}\label{eq:22}
\omega_n^0 = 0.
\end{equation}
But by (21)
$$
d\omega_n^0 = \omega_n^i \wedge (\mu_{ij} \omega_n^j 
+ \mu_i \omega_n^0).
$$
Hence the screen distribution  is integrable if and only if the tensor $\mu_{ij}$  is symmetric. Thus we arrived 
at the following result:

\begin{theorem} The fibration of normalizing subspaces $\zeta$ 
defines a screen distribution $\Delta$ of $(n-1)$-dimensional 
elements on a lightlike hypersurface $U^n$. This  distribution 
is integrable if and only if the tensor $\mu_{ij}$ 
defined by equation $(21)$ is symmetric.
\end{theorem}

Note that the configurations similar to that described in Theorem 
5 occurred in the works of the Moscow geometers published in the 
1950s. They were called the {\em one-side stratifiable pairs of 
ruled surfaces} (see \cite{F56}, \S 30 or \cite{AG93}, p. 187).

\section{An Invariant Normalization of 
Lightlike  Hypersurfaces of the de Sitter Space}

In \cite{A52} (see also \cite{AG96}, Ch. 2) 
an invariant normalization of a hypersurfaces $V^{n-1}$ 
of the conformal space $C^n$ was constructed. By Theorem 3, this normalization can be interpreted in terms of the geometry of 
the de Sitter space $S_1^{n+1}$.

Taking exterior derivative of equations (10) defining 
the hypersurface $V^{n-1}$ in the conformal space $C^n$, 
we obtain
$$
 \omega_i^n \wedge  \omega_0^i = 0,
$$
from which by linear independence of the 1-forms 
$\omega_0^i$ on $V^{n-1}$ and Cartan's lemma we find that 
\begin{equation}\label{eq:23}
\omega_i^n = \lambda_{ij} \omega_0^j, \;\; \lambda_{ij} 
= \lambda_{ji}.
\end{equation}
Here and in what follows we retain the notations used 
in the study of the geometry of hypersurfaces  
of the conformal space $C^n$ in the book \cite{AG96}. 

It is not difficult to find relations between the 
coefficients $\nu_{ij}$ in formulas (13) and 
 $\lambda_{ij}$ in formulas (23). Substituting the values 
of the forms $\omega_i^n$ and $\omega_0^j$ from (5) into 
(23), we find that 
$$
- g_{ij} \omega_n^j = \lambda_{ij} g^{jk} \omega_k^{n+1}.
$$
Solving these equations for $\omega_k^{n+1}$, we obtain 
$$
 \omega_i^{n+1} = - g_{ik} \widetilde{\lambda}^{kl} g_{lj}  \omega_n^j, 
$$
where $(\widetilde{\lambda}^{kl})$ is the inverse matrix of the 
matrix $(\lambda_{ij})$. Comparing these equations with equations (13), we obtain
\begin{equation}\label{eq:24}
\nu_{ij} = - g_{ik} \widetilde{\lambda}^{kl} g_{lj}.
\end{equation}
Of course, in this computation we assumed that the 
matrix $(\lambda_{ij})$ is nondegenerate.

Let us clarify the geometric meaning of the vanishing of 
$\det (\lambda_{ij})$. To this end, we make an admissible 
transformation of the moving frame associated with a point 
of a lightlike hypersurface $U^n$ by setting 
\begin{equation}\label{eq:25}
\widehat{A}_n = A_n + s A_0.
\end{equation}
The point $\widehat{A}_n$ as the point $A_n$ lies on the 
rectilinear generator $A_n A_0$. Differentiating this point and 
applying formulas (8) and (11), we obtain
\begin{equation}\label{eq:26}
d\widehat{A}_n = (ds + s\omega_0^0 + \omega_n^0) A_0 + 
(\omega_n^i + s \omega_0^i) A_i.
\end{equation}
It follows that in the new frame the forms $\omega_n^i$ 
become 
$$
\widehat{\omega}_n^i = \omega_n^i + s \omega_0^i.
$$
By (5) and (23), it follows that 
$$
\widehat{\omega}_n^i = -g^{ik} (\lambda_{kj} - s g_{kj}) \omega_0^j.
$$
This implies that in the new frame the quantities $\lambda_{ij}$ 
become
\begin{equation}\label{eq:27}
\widehat{\lambda}_{ij} = \lambda_{ij} - s g_{ij}.
\end{equation}
Consider also the matrix $(\widehat{\lambda}^i_j) 
= (g^{ik} \widehat{\lambda}_{kj})$. Since $g_{ij}$ is 
a nondegenerate tensor, the matrices  $(\widehat{\lambda}^i_j)$ 
and $(\widehat{\lambda}_{ij})$ have the same rank 
$\rho \leq n - 1$. 

From equation (26) it follows that 
$$
d\widehat{A}_n = (ds + s\omega_0^0 + \omega_n^0) A_0  
- \widehat{\lambda}^i_j A_i \omega^j_0.
$$
The differential $d\widehat{A}_n$ is the differential of the 
mapping $f: M^{n-1} \times l \rightarrow P^{n+1}$ 
which was considered in Theorem 1. The linearly independent 
forms $\omega_0^i$ are basis forms on 
the manifold $M^{n-1}$, and the form 
$\widehat{\omega}_n^0 = ds + s\omega_0^0 + \omega_n^0$ containing 
a nonhomogeneous parameter $s$ of the projective line $l$ is 
a basis form on this line. Thus the matrix 
\begin{equation}\label{eq:28}
\pmatrix{1 & 0 \cr
         0 & \widehat{\lambda}^i_j \cr}
\end{equation}
is the Jacobi matrix of this mapping.
Hence the tangent subspace to 
the  hypersurface $U^n$ at the point $\widehat{A}_n$ 
is determined by the points $\widehat{A}_n, A_0$, and 
$\widehat{\lambda}^i_j A_i$. At the points,  at which 
the rank $\rho$ of the 
matrix $(\widehat{\lambda}^i_j)$ is equal to $n-1, \rho = n - 1$, 
the tangent subspace to the  hypersurface $U^n$ has dimension $n$, and such points are {\em regular points} of the 
hypersurface. 
The points,  at which the rank $\rho$ of the 
matrix $(\widehat{\lambda}^i_j)$ is reduced, 
are {\em singular points} of the  hypersurface $U^n$. The 
coordinates 
of singular points are defined by the condition 
$\det (\widehat{\lambda}^i_j) = 0$ which by (27) is equivalent to 
the equation 
\begin{equation}\label{eq:29}
\det (\lambda_{ij} - s g_{ij}) = 0,
\end{equation}
the {\em characteristic equation} of the matrix $(\lambda_{ij})$ with respect to the tensor $g_{ij}$. The degree of this equation 
is equal to $n - 1$.

In particular, if  $A_n$ is a regular point 
of the  hypersurface $U^n$, then the matrix $(\lambda_{ij})$ 
is nondegenerate, and equation (24) holds. On the other hand, 
if  $A_n$ is a singular point 
of  $U^n$, then  equation (24) is meaningless. 

Since the matrix $(\lambda_{ij})$ is symmetric and 
the matrix $(g_{ij})$ defines a positive definite form 
of rank $n - 1$, equation (29) has 
$n - 1$ real roots if each root is counted as many times 
as its multiplicity. Thus on a rectilinear generator $A_n A_0$ of 
a lightlike   hypersurface $U^n$ there are $n-1$ real singular points. 

By Vieta's theorem, the sum of the roots of equation (29) is 
equal to the coefficient in $s^{n-2}$, and this coefficient 
is $\lambda_{ij} g^{ij}$. Consider the quantity
\begin{equation}\label{eq:30}
\lambda = \frac{1}{n-1} \lambda_{ij} g^{ij},
\end{equation}
which is the arithmetic mean of the roots of equation (29). 
This quantity $\lambda$ allows us to construct  new 
quantities
\begin{equation}\label{eq:31}
a_{ij} = \lambda_{ij} - \lambda g_{ij}.
\end{equation}
It is easy to check that the quantities $a_{ij}$ 
do not depend on the location of the point $A_n$ on the 
straight line $A_n A_0$, that is, $a_{ij}$ is invariant with 
respect to the transformation of the moving frame defined by 
equation (25). Thus the quantities $a_{ij}$ form a tensor 
on the  hypersurface $U^n$ defined in its 
second-order neighborhood. This tensor satisfies the condition 
\begin{equation}\label{eq:32}
a_{ij} g^{ij} = 0,
\end{equation}
that is, it is apolar to the tensor $g_{ij}$. 

On the straight line $A_n A_0$ we consider a point
\begin{equation}\label{eq:33}
C = A_n + \lambda A_0.
\end{equation}
It is not difficult to check that this point remains also 
fixed when the point $A_n$ moves along the  straight line $A_n 
A_0$. Since $\lambda$ is the arithmetic mean of the roots of 
equation (29) defining singular points on the  straight line $A_n A_0$, the point $C$ is the {\em harmonic pole} 
(see \cite{C50}) of the point $A_0$ with respect to these 
singular points. In particular, for $n = 3$, the point $C$ is the 
fourth harmonic point to the point $A_0$ with respect to two 
singular points of the rectilinear generator $A_3 A_0$ of the lightlike  hypersurface $U^3$ of the de Sitter space $S_1^4$. 

In the conformal theory of hypersurfaces, to 
the point $C$ there corresponds a hypersphere which is tangent 
to the hypersurface at the point $A_0$. This hypersphere
 is called the {\em central tangent hypersphere} 
(see \cite{AG96}, pp. 41--42). Since 
\begin{equation}\label{eq:34}
(d^2 A_0, C) = a_{ij} \omega_0^i \omega_0^j, 
\end{equation}
the cone 
$$
a_{ij} \omega_0^i \omega_0^j = 0
$$
with vertex at the point $A_0$ belonging to the tangent 
subspace $T_{A_0} (V^{n-1})$ contains the directions along which 
the central  hypersphere has a second-order tangency with 
the hypersurface $V^{n-1}$ at the point $A_0$. From the apolarity 
condition (33) it follows that it is possible to inscribe an orthogonal 
$(n-1)$-hedron with vertex at $A_0$ into the cone defined by equation (34). 

Now we can construct an invariant normalization of 
a lightlike hypersurface $U^n$ of the de Sitter space 
$S_1^{n+1}$. 
To this end, first we repeat some computations from Ch. 2 
of \cite{AG96}. 

Taking exterior derivatives of equations (23) 
and applying Cartan's lemma, we  obtain the equations
\begin{equation}\label{eq:35}
\nabla \lambda_{ij} +  \lambda_{ij} \omega_0^0 + g_{ij} \omega_n^0 =  \lambda_{ijk} \omega_0^k, 
\end{equation}
where 
$$
\nabla \lambda_{ij} = d \lambda_{ij} - \lambda_{ik} \omega^k_j 
-  \lambda_{kj} \omega^k_i, 
$$
and the quantities $ \lambda_{ijk}$ are symmetric 
with respect to all three indices. 
Equations (35) confirm one more time that the 
quantities $ \lambda_{ij} $ do not form a tensor 
and depend on a location of the point $A_n$ on the straight line 
$A_n A_0$. This dependence is described by a closed form 
relation (27). From formulas (35) it follows that the quantity 
$\lambda$ defined by equations (30) satisfy the differential 
equation
\begin{equation}\label{eq:36}
d \lambda  +  \lambda \omega_0^0 +  \omega_n^0 =  \lambda_k \omega_0^k, 
\end{equation}
where 
$$
\lambda_k = \frac{1}{n-1} g^{ij} \lambda_{ijk}
$$
(see formulas (2.1.36) and (2.1.37) in the book 
\cite{AG96}). The quantities $\lambda_k$ as well as 
the quantities $\lambda_{ijk}$ are determined by a 
third-order neighborhood of a generator $A_0 A_n$ of 
a lightlike hypersurface $U^n \subset S^{n+1}_1$.

The point $C$ lying on the rectilinear generator $A_n A_0$ 
of the  hypersurface $U^n$ describes a submanifold $W \subset U^n$ when  $A_n A_0$ moves. Let us find the tangent 
subspace to $U^n$ at the point $C$. Differentiating equation (33) 
and applying formulas (8) and (11), we obtain
$$
dC = (d \lambda + \lambda \omega_0^0 +  \omega_n^0) A_0 
+ (\omega_n^i +  \lambda \omega_0^i) A_i.
$$
By (5), (23), (30), and (36), it follows that 
 \begin{equation}\label{eq:37}
d C = (\lambda_i A_0  - g^{jk} a_{ki} A_j) \omega_0^i.
   \end{equation}
Define the affinor
 \begin{equation}\label{eq:38}
a^i_j = g^{ik} a_{kj}, 
 \end{equation}
whose rank coincides with the rank of the tensor $a_{ij}$. 
Then equation (37) takes the form
$$
d C = (\lambda_i A_0  -  a_i^j A_j) \omega_0^i.  
$$
The points
 \begin{equation}\label{eq:39}
C_i = \lambda_i A_0  -  a_i^j A_j
 \end{equation}
together with the point $C$ define the tangent subspace 
to the submanifold $W$ described by the point $C$ on 
 the  hypersurface $U^n$. 

If the point $C$ is a regular point of 
the rectilinear generator $A_n A_0$ of the  hypersurface $U^n$, 
then  the rank of the tensor $a_{ij}$ defined by 
equations (30) as well as the rank of the affinor $a_j^i$ is 
equal to $n-1$. As a result, the points $C_i$ are linearly 
independent and together with the point $C$ define 
the $(n-1)$-dimensional tangent subspace $T_C (W)$, and the submanifold $W$ itself has dimension $n-1, \dim W = n-1$.

The points $C_i$ also belong to the tangent subspace 
$T_{A_0} (V^{n-1})$ and define  the $(n-2)$-dimensional 
subspace $\zeta = T_{A_0} (V^{n-1}) \cap T_C (W)$ in it. This 
subspace 
is a normalizing subspace. Since such a normalizing subspace is 
defined in each tangent subspace $T_{A_0} (V^{n-1})$ of the 
hypersurface $V^{n-1} \subset Q^n$, there arises the fibration of 
these subspaces which by Theorem 4 defines an invariant affine 
connection on the lightlike hypersurface $U^n$. 
The subspace $T_C (W)$ is determined by a 
third-order neighborhood of a generator $A_0 A_n$ of 
the  hypersurface $U^n$.

 Thus we proved the following result:

\begin{theorem} If the tensor $a_{ij}$ defined by formula $(40)$ 
 on a lightlike hypersurface $U^n \subset S^{n+1}_1$ 
is nondegenerate, 
then it is possible to construct the invariant normalization 
of $U^n$ by means of the $(n-2)$-dimensional subspaces
$$
\zeta = C_1 \wedge C_2 \wedge \ldots \wedge C_{n-1}.
$$
This normalization induces on $U^n$ an invariant affine 
connection intrinsically connected with the geometry of this 
hypersurface. The normalization as well as the induced affine 
connection are determined by a 
third-order neighborhood of a generator $A_0 A_n$ of 
the  hypersurface $U^n$.
\end{theorem}

Theorem 5 implies that the invariant normalization 
we have constructed defines on $U^n$ an invariant screen 
distribution $\Delta$ which is also 
intrinsically connected with the geometry of the 
hypersurface $U^n$; here $\Delta_x = x \wedge \xi, \; x \in A_n A_0$.

Note that for the hypersurface $V^{n-1}$ of the conformal space $C^n$ 
a similar  invariant normalization was constructed as far 
back as 1953 (see \cite{A52} and also \cite{AG96}, Ch. 2). 
In the present paper we gave a new geometric meaning of 
this invariant normalization. 

\section{Singular Points of Lightlike  Hypersurfaces of 
the de Sitter Space}

As we indicated in Section 4, the points 
\begin{equation}\label{eq:40}
z =  A_n  + s A_0
\end{equation}
of the rectilinear generator $A_n A_0$ of the lightlike 
hypersurface $U^n$ are singular if their nonhomogeneous 
coordinate $s$ satisfies the equation 
\begin{equation}\label{eq:41}
\det (\lambda_{ij} - s g_{ij}) = 0,
\end{equation}
(In these points the Jacobian of the mapping $f: M^{n-1} \times l \rightarrow P^{n+1}$,  
which is equal to the determinant of the matrix (28), vanishes.)
We will investigate in more detail the structure of 
a lightlike hypersurface $U^n$ in a neighborhood of its singular 
point. 

Equation (41) is the characteristic equation of 
the matrix $(\lambda_{ij})$ with respect to the tensor 
$(g_{ij})$. The degree of this equation is $n-1$, and 
since the matrix  $(\lambda_{ij})$  is symmetric and 
 the matrix $(g_{ij})$ is also symmetric and 
positive definite, then according 
to the well-known result of linear algebra, all  roots 
of this equation are real, and the matrices $(\lambda_{ij})$  
 and  $(g_{ij})$ can be simultaneously reduced 
to a diagonal form. 

 Denote the roots of the characteristic equation 
 by $s_h, h = 1, 2, 
\ldots , n - 1$, and denote the corresponding singular points 
of the rectilinear generator $A_n A_0$ by 
\begin{equation}\label{eq:42}
B_h =  A_n  + s_h A_0.
\end{equation}
These singular points are called {\em foci}
 of the rectilinear generator $A_n A_0$ of a lightlike 
hypersurface $U^n$. 

It is clear from (42) that the point $A_0$ is not a focus 
 of the rectilinear generator $A_n A_0$. This is explained 
by the fact that by our assumption $\mbox{{\rm rank}} \; 
(\nu_{ij}) = n - 1$, and by (14), on the hyperquadric $Q^n$ 
the point $A_0$ describes a hypersurface $V^{n-1}$ which is 
transversal to the straight lines $A_0 A_n$. 

In the conformal theory of hypersurfaces, to 
the singular points $B_h$ there correspond the 
tangent hyperspheres defining the principal directions at 
a point $A_0$ of the hypersurface $V^{n-1}$ of the conformal 
space $C^n$ (see \cite{AG96}, p. 56). 

We will construct a classification of singular points 
of a lightlike hypersurface $U^n$ of the space $S^{n+1}_1$. 
We will use some computations that we made while constructing 
a classification of canal hypersurfaces in \cite{AG98a}. 

Suppose first that $B_1 = A_n + s_1 A_0$ be a singular 
point defined by a simple 
root $s_1$ of characteristic equation (41), 
$s_1 \neq s_h, h = 2, \ldots , n - 1$. For this singular point 
we have 
\begin{equation}\label{eq:43}
dB_1 = (ds_1 + s_1 \omega_0^0 + \omega_n^0) A_0 - 
\widehat{\lambda}_j^i \omega_0^j A_i,
\end{equation}
where 
\begin{equation}\label{eq:44}
\widehat{\lambda}_j^i = g^{ik} (\lambda_{kj} - s_1 g_{kj}) 
\end{equation}
is a degenerate symmetric affinor having a single null 
eigenvalue. The matrix of this affinor can be reduced 
to a quasidiagonal form
\begin{equation}\label{eq:45}
(\widehat{\lambda}_j^i) = \pmatrix{0 & 0 \cr 
             0 & \widehat{\lambda}_q^p \cr}, 
\end{equation}
where $p, q = 2, \ldots , n - 1$, and $(\widehat{\lambda}_q^p)$ 
is a nondegenerate symmetric affinor. The matrices 
$(g_{ij})$ and $(\lambda_{ij} - s_1 g_{ij})$ are reduced 
to the forms
$$
\pmatrix{1 & 0 \cr 
         0 & g_{pq}\cr} \;\; \mbox{{\rm and}} \;\;
\pmatrix{0 & 0 \cr 
             0 & \widehat{\lambda}_{pq}\cr}, 
$$
where $(\widehat{\lambda}_{pq}) = (\lambda_{pq} - s_1 
g_{pq})$ is a nondegenerate symmetric matrix. 

Since the point $B_1$ is defined invariantly on the 
generator $A_n A_0$, then it will be fixed if $\omega_0^i = 0$. 
Thus it follows from (43) that
\begin{equation}\label{eq:46}
ds_1 + s_1 \omega_0^0 + \omega_n^0 = s_{1i} \omega^i;
\end{equation}
 here and in what follows $\omega^i = \omega_0^i$. 
By (45) and (46) relation (43) takes the form 
\begin{equation}\label{eq:47}
dB_1 = s_{11} \omega^1 A_0 + (s_{1p} A_0 
- \widehat{\lambda}_p^q A_q) 
 \omega^p.
\end{equation}
Here the points $C_p = s_{1p} A_0 - \widehat{\lambda}_p^q A_q$ 
are linearly independent and belong to the tangent subspace 
$T_{A_0} (V^{n-1})$. 

Consider the submanifold ${\cal F}_1$ described by the singular point 
$B_1$ in the space $S^{n+1}_1$. This submanifold is called 
the {\em focal manifold} of the hypersurface $U^n$. 
Relation (47) shows that two cases are possible:

\begin{description}
\item[1)] $s_{11} \neq 0$. In this case  the submanifold 
${\cal F}_1$ 
is of dimension $n - 1$, and its tangent subspace at the point 
$B_1$ is determined by the points $B_1, A_0$, and $C_p$. This 
subspace contains the straight line $A_n A_0$, intersects 
the hyperquadric $Q^n$, and thus it,  as well as  the submanifold 
${\cal F}_1$ itself, is timelike. For $\omega^p = 0$, 
the point $B_1$ describes a curve $\gamma$ on the submanifold 
${\cal F}_1$  
which is tangent to the straight line $B_1 A_0$ coinciding with 
the generator $A_n A_0$ of the hypersurface $U^n$. 
The curve $\gamma$ is an isotropic curve of the de Sitter 
space $S^{n+1}_1$. Thus on 
${\cal F}_1$ there arises a fibre bundle of focal lines. The 
hypersurface $U^n$ is foliated into an $(n-2)$-parameter 
family of torses for which these lines are edges of regressions. 
The points $B_1$ are singular points of a kind which is called 
a {\em fold}. 

If the characteristic equation (41) has distinct roots, then 
an isotropic rectilinear generator $l$ of a lightlike hypersurface 
$U^n$ carries $n - 1$ distinct foci $B_h, h = 1, \ldots , n-1$. 
If for each of these foci the condition of type $s_{11} \neq 0$ 
holds, 
then each of them  describes a focal submanifold ${\cal F}_h$,  
carrying conjugate net.  Curves of one family of this net are 
tangent 
to the straight lines $l$, and this family is isotropic. 
On the hypersurface $V^{n-1}$ of the space $C^n = Q^n$ described by 
the point $A_0$, to these conjugate nets there correspond the 
net of curvature lines. 

\item[2)] $s_{11} = 0$. In this case  relation (47) takes 
the form
\begin{equation}\label{eq:48}
dB_1 =  (s_{1p} A_0 - \widehat{\lambda}_p^q A_q)  \omega^p,
\end{equation}
and the focal submanifold ${\cal F}_1$ 
is of dimension $n - 2$. Its tangent subspace at the point 
$B_1$ is determined by the points $B_1$ and $C_p$. 
An arbitrary point $z$ of this subspace can be written in the 
form 
$$
z = z^n B_1 + z^p C_p = z^n (A_n + s_1 A_0) + z^p 
(s_{1p} A_0 - \widehat{\lambda}_p^qA_q).
$$
Substituting the coordinates of this point into relation (3), 
we find that 
$$
(z, z) = g_{rs} \widehat{\lambda}_p^r   
\widehat{\lambda}_q^s z^p z^q + (z^n)^2 > 0.
$$
It follows that the tangent subspace 
$T_{B_1} ({\cal F}_1)$ does not have common points with 
the hyperquadric $Q^n$, that is, it is spacelike. Since this 
takes place for any point $B_1 \in {\cal F}_1$, the focal 
submanifold ${\cal F}_1$ is spacelike.

For $\omega^p = 0$, the point $B_1$ is fixed. The subspace 
$T_{B_1} ({\cal F}_1)$ will be fixed too. 
On the hyperquadric $Q^n$, the point $A_0$ describes 
a curve $q$ which is polar-conjugate to $T_{B_1} ({\cal F}_1)$. 
Since $\dim T_{B_1} ({\cal F}_1) = n - 2$, the curve $q$ is 
a conic, along which the two-dimensional plane 
polar-conjugate to the subspace $T_{B_1} ({\cal F}_1)$  with 
respect to the hyperquadric $Q^n$, intersects $Q^n$. 
Thus for $\omega^p = 0$, the rectilinear generator $A_n A_0$ of 
the hypersurface $U^n$ describes  a two-dimensional 
second-order cone  with vertex at the point $B_1$ 
and the directrix $q$. Hence in the case under consideration 
a lightlike hypersurface $U^n$ is foliated into an 
$(n-2)$-parameter family of second-order cones whose vertices 
describe the $(n-2)$-dimensional focal submanifold 
${\cal F}_1$, and 
the points $B_1$ are {\em conic} singular points of 
the hypersurface $U^n$. 

The hypersurface $V^{n-1}$ of the conformal space $C^n$ 
corresponding to such a lightlike hypersurface $U^n$ 
is a canal hypersurface which envelops an $(n-2)$-parameter 
family of hyperspheres. Such a hypersurface carries a family 
of cyclic generators which depends on the same number 
of parameters. Such hypersurfaces were investigated in detail 
in \cite{AG98a}. 
\end{description}

Further let $B_1$ be a  singular point of multiplicity $m$, where 
$m \geq 2$, of a rectilinear generator $A_n A_0$ 
of  a lightlike hypersurface $U^n$ of the space $S^{n+1}_1$ 
defined by an $m$-multiple root of characteristic equation 
(41). We will assume that 
\begin{equation}\label{eq:49}
 s_1 = s_2 = \ldots = s_m := s_0, s_0 \neq s_p, 
\end{equation}
and also assume that $a, b, c = 1, \ldots, m$ 
and $p, q, r = m+1, \ldots, n-1$. Then the matrices $(g_{ij})$ 
and $(\lambda_{ij})$ can be simultaneously  
reduced to  quasidiagonal forms
$$
\pmatrix{g_{ab} & 0 \cr 
             0 & g_{pq}\cr} \;\; \mbox{{\rm and}} \;\;
\pmatrix{s_0 g_{ab} & 0 \cr 
             0 & \lambda_{pq}\cr}.
$$
We also construct the matrix $(\widehat{\lambda}_{ij}) 
= (\lambda_{ij} - s_0 g_{ij})$. Then 
\begin{equation}\label{eq:50}
(\widehat{\lambda}_{ij}) = \pmatrix{0 & 0 \cr 
             0 & \widehat{\lambda}_{pq}\cr},
\end{equation}
where $\widehat{\lambda}_{pq} = \lambda_{pq} - s_0 g_{pq}$ 
is a nondegenerate matrix of order $n - m - 1$. 

By relations (50) and formulas (5) and (23) we have
\begin{equation}\label{eq:51}
\omega_a^n - s_0 \omega_a^{n+1} = 0,
\end{equation}
\begin{equation}\label{eq:60}
\omega_p^n - s_0 \omega_p^{n+1} =  \widehat{\lambda}_{pq} 
\omega^q.
\end{equation}
Taking exterior derivative of equation (51) and applying relation 
(52), we find that 
\begin{equation}\label{eq:53}
\widehat{\lambda}_{pq} \omega_a^p \wedge \omega^q 
+ g_{ab} \omega^b \wedge (ds_0 + s_0 \omega_0^0 
+ \omega_n^0) = 0.
\end{equation}
It follows that the 1-form $ds_0 + s_0 \omega_0^0 
+ \omega_n^0$ can be expressed in terms of the basis forms. 
We write these expressions in the form
\begin{equation}\label{eq:54}
ds_0 + s_0 \omega_0^0 + \omega_n^0 = s_{0c} \omega^c + 
s_{0q} \omega^q.
\end{equation}
Substituting this decomposition into equation (53), we find that 
\begin{equation}\label{eq:55}
(\widehat{\lambda}_{pq} \omega_a^p  
+ g_{ab} s_{0q} \omega^b) \wedge \omega^q + g_{ab}  s_{0c} 
\omega^b \wedge \omega^c  = 0.
\end{equation}
The terms in the left hand side of (55) do not have similar 
terms. Hence both terms are equal to 0. Equating to 0 
the coefficients of the summands of the second term, 
we find that 
\begin{equation}\label{eq:56}
 g_{ab}  s_{0c} =  g_{ac}  s_{0b}.
\end{equation}
Contracting this equation with the matrix $(g^{ab})$ 
which is the inverse matrix of the matrix $(g_{ab})$, we obtain
$$
ms_{0c} = s_{0c}.
$$
Since $m \geq 2$, it follows that 
$$
s_{0c} = 0,
$$
and  relation (54) takes the form 
\begin{equation}\label{eq:57}
ds_0 + s_0 \omega_0^0 + \omega_n^0 = s_{0p} \omega^p.
\end{equation}

For the  singular point of multiplicity $m$ of the generator 
$A_n A_0$ in question the equation (43) can be written in the 
form
$$
dB_1 = (ds_0 + s_0 \omega_0^0 + \omega_n^0) A_0 - 
\widehat{\lambda}_q^p \omega_0^q A_p.
$$
Substituting decomposition (57) in the last equation, 
we find that 
\begin{equation}\label{eq:58}
dB_1 = (s_{0p} A_0 - \widehat{\lambda}_p^q A_q)\omega_0^p.
\end{equation}
This relation is similar to equation (48) with the only 
difference that in (48) we had $p, q = 2, \ldots, n - 1$, and 
in (58) we have $p, q = m+1, \ldots, n - 1$. 
Thus the point $B_1$ describes now a spacelike focal 
manifold ${\cal F}_1$ of dimension $n-m-1$. For $\omega_0^p = 0$, 
the point $B_1$ is fixed, and the point $A_0$ describes an 
$m$-dimensional submanifold on the hyperquadric $Q^n$ which 
is a cross-section of $Q^n$ by an $(m+1)$-dimensional subspace 
that is polar-conjugate to the $(n-m-1)$-dimensional subspace 
tangent to the submanifold ${\cal F}_1$. 

The point $B_1$ is a  conic singular point of multiplicity $m$ of 
a lightlike hypersurface $U^n$, and this hypersurface is foliated 
into an $(n-m-1)$-parameter family of $(m+1)$-dimensional 
second-order cones circumscribed about 
the hyperquadric $Q^n$. The  hypersurface $V^{n-1}$ of the 
conformal space $C^n$ that corresponds to such a  hypersurface 
$U^n$ is an $m$-canal hypersurface (i.e.,  
the envelope of an $(n-m-1)$-parameter family of 
hyperspheres), and it carries an $m$-dimensional 
spherical generators.

Note also an extreme case when the rectilinear generator 
$A_n A_0$ of a lightlike hypersurface $U^n$ carries a single 
 singular point of multiplicity $n-1$. 
As follows from our consideration 
of the cases $m \geq 2$, this singular point is fixed, and 
the hypersurface $U^n$ become a second-order hypercone with 
vertex at this singular point which is circumscribed about 
 the hyperquadric $Q^n$. This hypercone is the isotropic cone 
of the space $S^{n+1}_1$. The  hypersurface $V^{n-1}$ of the 
conformal space $C^n$ that corresponds to such a  hypersurface 
$U^n$ is a hypersphere of the space $C^n$.
 
The following theorem combines the results of this section:

\begin{theorem}
 A lightlike hypersurface $U^n$ of maximal rank $r=n-1$ of the de 
Sitter space $S^{n+1}_1$ possesses $n-1$ real singular points 
on each of its rectilinear generators if each of these singular 
points  is counted as many times as its multiplicity. The simple 
singular points can be of two kinds: a fold and conic. In the 
first case the hypersurface $U^n$ is foliated into an $(n-2)$-parameter family of torses, and in the second case it is foliated 
into an $(n-2)$-parameter family of second-order cones. 
The vertices of these cones 
describe the $(n-2)$-dimensional spacelike submanifold 
in the space $S^{n+1}_1$. All multiple singular points 
of a hypersurface $U^n$ are conic. If a rectilinear 
generator of a hypersurface $U^n$ carries a singular point 
of multiplicity $m$, $2 \leq m \leq n -1$, 
then the hypersurface $U^n$ is foliated into an $(n-m-1)$-parameter family of $(m+1)$-dimensional second-order cones. 
The vertices of these cones 
describe the $(n-m-1)$-dimensional spacelike submanifold 
in the space $S^{n+1}_1$.
The hypersurface $V^{n-1}$ of the conformal space $C^n$ 
corresponding to a lightlike hypersurface $U^n$ with 
 singular points of multiplicity $m$ 
is a canal hypersurface which envelops an $(n-m-1)$-parameter 
family of hyperspheres and has $m$-dimensional spherical 
generators. 
\end{theorem}

Since lightlike hypersurfaces $U^n$ of 
the de Sitter space $S^{n+1}_1$ represent a light flux 
(see Section {\bf 2}), its focal submanifolds have the following 
physical meaning. If one of them is a lighting submanifold, then 
others will be manifolds of concentration of a light flux. 
Intensity of concentration depends on multiplicity of a focus 
 describing this submanifold. 

In the extreme case when an isotropic rectilinear generator 
$l = A_n A_0$ of a hypersurface $U^n$ carries one $(n-1)$-multiple 
focus, the hypersurfaces $U^n$ degenerates into the light cone 
generated by a point source of light. This cone represents a radiating 
light flux.

If each isotropic generator $l \subset U^n$ carries two foci 
$B_1$ and $B_2$ of multiplicities $m_1$ and $m_2, \,m_1 + m_2 = n - 1, \,
m_1 > 1, \, m_2 > 1$, then these foci describe spacelike submanifolds 
${\cal F}_1$ and ${\cal F}_2$ of dimension $n - m_1 - 1$ and 
$n - m_2 - 1$, respectively. If one of these submanifolds is a lighting 
 submanifold, then on the second one a light flux is concentrated.

\noindent
{\sf             
Department of Mathematics,            
Jerusalem College of Technology - Mahon Lev,  
21 Havaad Haleumi St., POB 16031,  Jerusalem  91160, Israel       
}

\noindent
{\em E-mail address}: akivis@math.jct.ac.il

\vspace*{3mm}

\noindent
{\sf 
 Department of Mathematics,   New Jersey Institute of Technology,  
 University Heights, Newark, NJ 07102
}

\noindent
{\em E-mail address}: vlgold@numerics.njit.edu
\end{document}